\documentclass[12pt]{article}
\usepackage{amssymb, amsfonts, amssymb, amsthm, amscd}

\newtheorem{theorem}{{\bf Theorem}}[section]

\newtheorem{proposition}{{\bf Proposition}}[section]
\newtheorem{definition}{{\bf Definition}}[section]

\newtheorem{remark}{{\bf Remark}}[section]

\input tcilatex
\begin{document}

\title{Constrained BSDE and Viscosity Solutions of Variation Inequalities\footnote{This work is supported by the National
Basic Research Program of China (973 Program), No. 2007CB814902 and
No. 2007CB814906.}}
\date{August 23, 2007}
\author{Shige PENG$^{a,c}$\ \ \
Mingyu XU $^{a,c}$\footnote{Corresponding author, Email:
xumy@amss.ac.cn}
\\
%EndAName
{\small $^a$School of Mathematics and System Science, Shandong University, }%
\\
{\small 250100, Jinan, China}\\
{\small $^b$Institute of Applied Mathematics, Academy of Mathematics
and Systems
Science,}\\
{\small Chinese Academy of Sciences, Beijing, 100080, China.}\\
{\small $^c$Department of Financial Mathematics and Control science,
School of
Mathematical Science,} \\
{\small Fudan University, Shanghai, 200433, China.}} \maketitle

\textbf{Abstract.} In this paper, we study the relation between the
smallest $g$-supersolution of constraint backward stochastic
differential equation and viscosity solution of constraint
semilineare parabolic PDE, i.e. variation inequalities. And we get
an existence result of variation inequalities via constraint BSDE,
and prove a uniqueness result under certain condition.
 \\[.5cm]

\section{Introduction}

El Karoui, Kapoudjian, Pardoux, Peng and Quenez (1997) studied the problem
of BSDE (backward stochastic differential equation) with reflection, which
is, a standard BSDE with an additional continuous, increasing process in
this equation to keep the solution above a certain given continuous boundary
process. This increasing process must be chosen in certain minimal way, i.e.
an integral condition, called Skorohod reflecting condition (cf. \cite
{Skoro1965}), is satisfied. It was proved in this paper that the solution of
the reflected BSDE associated to a terminal condition $\xi $, a coefficient $%
g$ and a lower reflecting obstacle $L$, is the smallest supersolution of
BSDE with same parameter $(\xi ,g)$, which dominates the given boundary
process $L$. Then in same paper, they give a probabilistic interpretation of
viscosity solution of variation inequality by the solution of reflected
BSDEs.

An important application of the constrained BSDE is the pricing of
contingent claims with constraint of protfolios, i.e. portfolios of an asset
is constrained in a given subset. In this case the solution $(y,z)$ of the
corresponding reflected BSDE must remain in this subset. This problem was
studied by Karaztas and Kou (cf. \cite{KK}), then by \cite{EQ95} and \cite
{CKS}.

The most general case of the constraint $\Phi $, which is discribed by a
Lipschitz continuous function, is first studied in \cite{Peng99}. Author
proved that under the Lipschitz condition of coefficient $g$, the smallest
supersolution of BSDE with coefficient $g$ and constraint $\Phi \geq 0$
exists, if there exists a special solutin of this constraint problem.

The main conditions of our paper is same as \cite{Peng99}: $g$ is a
Lipschitz function and the constraint $\Phi (\omega ,t,x,y,z)$, $t\in [0,T]$
is a Lipschitz continuous function. In this paper we study the relation
between the smallest $g$-supersolution and viscosity solution of constraint
semilineare parabolic PDE, i.e. variation inequalities. And we get an
existence result of variation inequalities via constraint BSDE, and prove a
uniqueness result under certain condition.

\section{Preliminaries and Constraint BSDEs}

Let $(\Omega ,\mathcal{F},P)$ be a probability space, and $%
B=(B_{1},B_{2},\cdots ,B_{d})^{T}$ be a $d$-dimensional Brownian motion
defined on $[0,\infty )$. We denote by $\{\mathcal{F}_{t};0\leq t<\infty \}$
the natural filtration generated by this Brownian motion $B:$%
\[
\mathcal{F}_{t}=\sigma \{\{B_{s};0\leq s\leq t\}\cup \mathcal{N\}},
\]
where $\mathcal{N}$ is the collection of all $P-$null sets of $\mathcal{F}$.
The Euclidean norm of an element $x\in \mathbb{R}^{m}$ is denoted by $|x|$.
We also need the following notations for $p\in [1,\infty )$:

\begin{itemize}
\item  $\mathbf{L}^{p}(\mathcal{F}_{t};\mathbb{R}^{m}):=$\{$\mathbb{R}^{m}$%
-valued $\mathcal{F}_{t}$--measurable random variables $X$ s.t. $%
E[|X|^{p}]<\infty $\};

\item  $\mathbf{L}_{\mathcal{F}}^{p}(0,t;\mathbb{R}^{m}):=$\{$\mathbb{R}^{m}$%
--valued and $\mathcal{F}_{t}$--progressively measurable processes $\varphi $
defined on $[0,t]$, s.t. $E\int_{0}^{t}|\varphi _{s}|^{p}ds<\infty $\};

\item  $\mathbf{D}_{\mathcal{F}}^{p}(0,t;\mathbb{R}^{m}):=$\{$\mathbb{R}^{m}$%
--valued and RCLL $\mathcal{F}_{t}$--progressively measurable processes $%
\varphi $ defined on $[0,t]$, s.t. $E[\sup_{0\leq s\leq t}|\varphi
_{s}|^{p}]<\infty $\};

\item  $\mathbf{A}_{\mathcal{F}}^{p}(0,t):=$\{increasing processes $A$ in $%
\mathbf{D}_{\mathcal{F}}^{p}(0,t;\mathbb{R})$ with $A(0)=0$\}.
\end{itemize}

\noindent When $m=1$, they are denoted by $\mathbf{L}^{p}(\mathcal{F}_{t})$,
$\mathbf{L}_{\mathcal{F}}^{p}(0,t)$ and $\mathbf{D}_{\mathcal{F}}^{p}(0,t)$,
respectively. We are mainly interested in the case $p=2$. In this paper, we
consider BSDE on the interval $[0,T]$, with a fixed $T>0$.

We put the BSDE with a constraint into Markovian framework. Consider the
following forward SDE,
\begin{eqnarray}
dX_{s}^{t,x} &=&b(s,X_{s}^{t,x})ds+\sigma (s,X_{s}^{t,x})dW_{s},\;t\leq
s\leq T,  \label{SDE1} \\
X_{t}^{t,x} &=&x.  \nonumber
\end{eqnarray}
where $b:[0,T]\times \mathbf{R}^{d}\rightarrow \mathbf{R}^{d}$, $\sigma
:[0,T]\times \mathbf{R}^{d}\rightarrow \mathbf{R}^{d\times d}$ are
continuous mappings, satisfying

\[
\begin{array}{rll}
\mbox{(i)} & b\mbox{ and }\sigma \mbox{ are continuous in }t &  \\
\mbox{(ii)} & \left| b(t,x)-b(t,x^{\prime })\right| +\left| \sigma
(t,x)-\sigma (t,x^{\prime })\right| \leq k(\left| x-x\right| ), & dP\times dt%
\mbox{ a.s. }
\end{array}
\]
for some $k>0$, and for all $x$, $x^{\prime }\in \mathbf{R}^{d}$. And for
each $(t,x)\in [0,T]\times \mathbf{R}^{d}$, $\{X_{s}^{t,x};t\leq s\leq T\}$
is denoted as the unique solution of SDE (\ref{SDE1}).

Let $g$ be a coefficient $g(t,x,y,z):[0,T]\times \mathbb{ R}^{d}\times %
\mathbb{R\times R}^{d}\rightarrow \mathbb{R}$, which satisfies the following
assumptions: there exists a constant $\mu >0$, $p\in \mathbf{N}$ such that,
for each $x$ in $\mathbb{R}^{d}$, $y,y^{\prime }\;$in $\mathbb{R}$ and $%
z,z^{\prime }$ in $\mathbb{R}^{d}$, we have
\begin{equation}
\begin{array}{rll}
\mbox{(i)} & \left| g(t,x,0,0)\right| \leq \mu (1+\left| x\right| ^{p}) &
\\
\mbox{(ii)} & \left| g(t,x,y,z)-g(t,x,y^{\prime },z^{\prime })\right| \leq
\mu (\left| y-y^{\prime }\right| +\left| z-z^{\prime }\right| ), & dP\times
dt\mbox{ a.s. }
\end{array}
\label{Lip}
\end{equation}

Our BSDE with a constraint is
\begin{eqnarray}
-dY_{s}^{t,x}
&=&g(s,X_{s}^{t,x},Y_{s}^{t,x},Z_{s}^{t,x}{}_{s})ds+dA_{s}^{t,x}-Z_{s}^{t,x}dW_{s},
\label{CBSDE} \\
Y_{T}^{t,x} &=&\Psi (X_{T}^{t,x}),\;\;\;\mbox{with }\Phi
(s,X_{s}^{t,x},Y_{s}^{t,x},Z_{s}^{t,x})\geq 0,d\mathbf{P}\times dt\mbox{%
-a.s..}  \nonumber
\end{eqnarray}
Here $\Psi :\mathbf{R}^{d}\rightarrow \mathbf{R}$, has at most polynomial
growth at infinity. $\Phi :[0,T]\times \mathbf{R}^{d}\times \mathbf{R\times R%
}^{d}\rightarrow \mathbf{R}$, which plays a role of constraint in this
paper, satisfying: there exists a constant $\mu _{2}>0$, such that, for each
$x\in \mathbb{R}^{d}$, $y$, $y^{\prime }\in \mathbb{R}$ and $z$, $z^{\prime
}\in \mathbb{R}^{d}$, we have
\begin{equation}
\begin{array}{lll}
\mbox{(i)} & \left| \Phi (t,x,0,0)\right| \leq \mu _{2}(1+\left| x\right|
^{p}) &  \\
\mbox{(ii)} & \left| \Phi (t,x,y,z)-\Phi (t,x,y^{\prime },z^{\prime
})\right| \leq \mu _{2}(\left| y-y^{\prime }\right| +\left| z-z^{\prime
}\right| ), & dP\times dt\mbox{ a.s. } \\
\mbox{(ii)} & y\rightarrow \Phi (t,x,\cdot ,z)\mbox{ and
}z\rightarrow \Phi (t,x,y,\cdot )\mbox{ are continuous.} &
\end{array}
\label{Constraint}
\end{equation}
The constraint $\Phi $ is an equivalent form of the constraint we have
discussed before, as \cite{EKPPQ}, \cite{Peng99} and \cite{PX06}.

\begin{definition}
The solution of (\ref{CBSDE}) is $(Y_{s}^{t,x},Z_{s}^{t,x},A_{s}^{t,x})_{t%
\leq s\leq T}$ defined as the smallest $g$--supersolution constrained by $%
\Phi \geq 0$, i.e. $Y^{t,x}\in \mathbf{D}_{\mathcal{F}}^{2}(t,T)$ and there
exist a predictable process $Z^{t,x}\in \mathbf{L}_{\mathcal{F}}^{2}(t,T;%
\mathbb{R}^{d})$ and an increasing RCLL process $A^{t,x}\in \mathbf{A}_{%
\mathcal{F}}^{2}(t,T)$ such that (\ref{CBSDE}) is satisfied and if there is
another process $Y^{t,x\prime }\in \mathbf{D}_{\mathcal{F}}^{2}(t,T)$, with $%
(Z^{t,x\prime },A^{t,x\prime })\in \mathbf{L}_{\mathcal{F}}^{2}(t,T;%
\mathbb{R}^{d})\times \mathbf{A}_{\mathcal{F}}^{2}(t,T)$, satisfying (\ref
{CBSDE}), then we have $Y_{s}^{t,x\prime }\geq Y_{s}^{t,x}$.
\end{definition}

The following theorem of the existence of the smallest solution was obtained
in \cite{Peng99}.

\begin{theorem}
\label{exist} Suppose that $\xi \in \mathbf{L}^{2}(\mathcal{F}_{T})$, the
function $g$ satisfies (\ref{Lip}) and the constraint $\Phi $ satisfies (\ref
{Constraint}). We assume that \textbf{(H) }there is one $g$--supersolution $%
y^{\prime }\in \mathbf{D}_{\mathcal{F}}^{2}(0,T)$, constrained by $\Phi \geq
0$:
\begin{eqnarray}
y_{t}^{\prime } &=&\xi +\int_{t}^{T}g(s,y_{s}^{\prime },z_{s}^{\prime
})ds+A_{T}^{\prime }-A_{t}^{\prime }-\int_{t}^{T}z_{s}^{\prime }dB_{s},
\label{no-empty1} \\
A^{\prime } &\in &\mathbf{A}_{\mathcal{F}}^{2}(0,T)\mbox{ ,\ }\Phi
(t,y_{t}^{\prime },z_{t}^{\prime })\geq 0,\;dP\times dt\mbox{ a.s. }
\nonumber
\end{eqnarray}
Then there exists the smallest $g$--supersolution $y\in \mathbf{D}_{\mathcal{%
F}}^{2}(0,T)$ constrained by $\Phi \geq 0$, with the terminal condition $%
y_{T}=\xi $, i.e. there exists a triple $(y_{t},z_{t},A_{t})\in \mathbf{D}_{%
\mathcal{F}}^{2}(t,T)\times \mathbf{L}_{\mathcal{F}}^{2}(t,T;\mathbb{R}%
^{d})\times \mathbf{A}_{\mathcal{F}}^{2}(t,T)$, such that
\begin{eqnarray*}
y_{t} &=&\xi
+\int_{t}^{T}g(s,y_{s},z_{s})ds+A_{T}-A_{t}-\int_{t}^{T}z_{s}dB_{s}, \\
A &\in &\mathbf{A}_{\mathcal{F}}^{2}(0,T)\mbox{ ,\ }\Phi (t,y_{t},z_{t})\geq
0,\;dP\times dt\mbox{ a.s. }
\end{eqnarray*}
Moreover, this smallest $g$--supersolution is the limit of a sequence of $%
g^{n}$--solutions with $g^{n}=g+n\Phi ^{-}$, where the convergence is in the
following sense:
\begin{eqnarray}
y_{t}^{n} &\nearrow &y_{t}\mbox{, with }\lim_{n\rightarrow \infty
}E[|y_{t}^{n}-y_{t}|^{2}]=0,\;\lim_{n\rightarrow \infty
}E\int_{0}^{T}|z_{t}-z_{t}^{n}|^{p}dt=0,p\in [1,2),\;  \label{approx1} \\
A_{t}^{n} &:&=\int_{0}^{t}(g+n\Phi ^{-})(s,y_{s}^{n},z_{s}^{n})ds\rightarrow
A_{t}\mbox{ weakly in }\mathbf{L}^{2}(\mathcal{F}_{t}),  \label{approx2}
\end{eqnarray}
where $z$ and $A$ are the corresponding martingale part and increasing part
of $y$, respectively.
\end{theorem}

And we recall an interesting proposition proved in \cite{PX06}.

\begin{proposition}
\label{gmphi}A process $y\in \mathbf{D}_{\mathcal{F}}^{2}(0,T)$ is
the smallest $g$-supersolution on $[0,T]$ constraint by $\Phi$ with
$y_T=\xi$, if and only if for all $m\geq 0$, it is a
$(g+m\Phi)$-supersolution on $[0,T]$ with $y_T=\xi$.
\end{proposition}
\section{Relation between BSDE with a constraint and PDE}

In the following, we assume that \textbf{(H)} holds, and denote the smallest
solution of (\ref{CBSDE}) by $(Y_{s}^{t,x},Z_{s}^{t,x},A_{s}^{t,x})_{t\leq
s\leq T}$. Define
\[
u(t,x):=Y_{t}^{t,x}.
\]
The variation inequality we concerned is
\[
\min \{-\partial _{t}u-F_{0}(t,x,u,Du,D^{2}u),\Phi (x,u,\sigma
^{T}(x)Du)\}=0,
\]
where $F_{0}(t,x,u,q,S):=\frac{1}{2}\sum_{i,j=1}^{n}[\sigma \sigma
^{T}]_{ij}(t,x)S_{ij}+\left\langle b(t,x),q\right\rangle +g(t,x,u,\sigma
^{T}(t,x)q)$. We study this problem by the following penalization approach:
for each $\alpha \geq 0$,
\begin{eqnarray*}
-dY_{s}^{t,x,\alpha } &=&(g+\alpha \Phi
^{-})(s,X_{s}^{t,x},Y_{s}^{t,x,\alpha },Z_{s}^{t,x,\alpha
})ds-Z_{s}^{t,x,\alpha }dW_{s}, \\
Y_{T}^{t,x,\alpha } &=&\Psi (X_{T}^{t,x}).\;\;\;
\end{eqnarray*}
Define
\begin{equation}
u_{\alpha }(t,x):=Y_{t}^{t,x,\alpha }.  \label{eq4.6}
\end{equation}
Then by theorem \ref{exist}, we have
\begin{equation}
u_{\alpha }(t,x)\nearrow u(t,x),\;(t,x)\in [0,T]\times
\mathbf{R}^{n},\mbox{ as }\alpha \rightarrow \infty .  \label{eq4.8}
\end{equation}

We introduce the following penalized PDE
\begin{equation}
\partial _{t}u+F_{\alpha }(t,x,u,Du,D^{2}u)=0,\;\;\forall \alpha =1,2,\cdots
,  \label{eq4.9}
\end{equation}
where $F_{\alpha }(t,x,u,q,S):=\frac{1}{2}\sum_{i,j=1}^{n}[\sigma \sigma
^{T}]_{ij}(t,x)S_{ij}+\left\langle b(t,x),q\right\rangle +(\alpha \Phi
^{-}+g)(t,x,u,\sigma ^{T}(t,x)q)$.

To introduce the definition of viscosity solution. First we need the notions
of parabolic superjet and subjet.

\begin{definition}
For\mbox{ a} function $u\in $LSC$([0,T]\times R^{n})$ (resp. USC$%
([0,T]\times R^{n})$), we define the parabolic superjet (resp. parabolic
subjet) of $u$ at $(t,x)$ by $\mathcal{P}^{2,+}u(t,x)$ (resp. $\mathcal{P}%
^{2,-}u(t,x)$), the set of triples $(p,q,X)\in \mathbf{R\times R}^{d}\mathbf{%
\times S}^{n}$, satisfying
\begin{eqnarray*}
u(s,y) &\leq &\mbox{(resp.}\geq \mbox{) }u(t,x)+p(s-t)+\left\langle
q,y-s\right\rangle +\frac{1}{2}\left\langle X(y-x),y-x\right\rangle  \\
&&+o\left( \left| s-t\right| +\left| y-x\right| ^{2}\right) .
\end{eqnarray*}
\end{definition}

Then we have

\begin{definition}
\mbox{ A} function $u\in $LSC$([0,T]\times R^{n})$ (resp.
USC$([0,T]\times
R^{n})$) is called a viscosity supersolution (resp. subsolution) of $%
\partial _{t}u+F_{\alpha }=0$ if for each $(t,x)\in (0,T)\times R^{n}$, for
any $(p,q,X)\in \mathcal{P}^{2,-}u(t,x)$(resp. $(p,q,X)\in \mathcal{P}%
^{2,+}u(t,x)$), we have
\[
p+F_{\alpha }(t,x,u(t,x),q,X)\leq 0,\;\mbox{(resp.\ }\geq 0\mbox{).}
\]
\end{definition}

The following result can be found in \cite{P91}.

\begin{proposition}
\label{u-alpha}\mbox{ }We assume (\ref{Lip}) and $\Psi $ has at most
polynomial growth at infinity. Then for each $\alpha =1,2,\cdots $,
the function $u_{\alpha }\in $C$([0,T]\times R^{n})$ defined by
(\ref{eq4.6}) is the viscosity solution of $\partial _{t}u_{\alpha
}+F_{\alpha }=0$.
\end{proposition}

Now we return to the variation inequality
\begin{equation}
\min \{-\partial _{t}u-F_{0}(t,x,u,Du,D^{2}u),\Phi (x,u,\sigma
^{T}(x)Du)\}=0.  \label{PDECons}
\end{equation}
The solution of this equation may be not continuous, so we need the
definition of discontinuous viscosity solution. For a given locally bounded
function $v$, we define its upper and lower semicontinuous envelope of $v$,
denoted as $v^{*}$ and $v_{*}$ respectively, where
\[
v^{*}(t,x)=\limsup_{t^{\prime }\rightarrow t,x^{\prime }\rightarrow
x}v(t^{\prime },x^{\prime }),v_{*}(t,x)=\liminf_{t^{\prime }\rightarrow
t,x^{\prime }\rightarrow x}v(t^{\prime },x^{\prime }).
\]
Then

\begin{definition}
\label{def4.2}(i)\mbox{ A} locally bounded function $u$ is called a
viscosity supersolution (\ref{PDECons}) if for each $(t,x)\in
(0,T)\times \mathbf{R}^{n}$, for any $(p,q,X)\in
\mathcal{P}^{2,-}u_{*}(t,x)$, then we have
\begin{equation}
\min \{-p-F_{0}(t,x,u_{*},q,X),\Phi (x,u_{*},\sigma ^{T}(x)q)\}\geq 0,
\label{VSCsuper}
\end{equation}
i.e. we have both $\Phi (x,u_{*},\sigma ^{T}(x)q)\geq 0$ and
\[
-p-F_{0}(t,x,u_{*},q,X)\geq 0.
\]
(ii) \mbox{A} locally bounded function $u$ is called a viscosity
subsolution of (\ref{PDECons}), if for each $(t,x)\in (0,T)\times
\mathbf{R}^{n}$, for any $(p,q,X)\in \mathcal{P}^{2,+}u^{*}(t,x)$,
then we have
\begin{equation}
\min \{-p-F_{0}(t,x,u^{*},q,X),\Phi (x,u^{*},\sigma ^{T}(x)q)\}\leq 0,
\label{VSCsub}
\end{equation}
i.e. for $(t,x)\in (0,T)\times R^{n}$ where $\Phi (x,u^{*},\sigma ^{T}(x)q)>0
$, we have
\[
-p-F_{0}(t,x,u^{*},q,X)\leq 0.
\]
(iii) \mbox{A} locally bounded function $u$ is called a viscosity
solution of (\ref{PDECons}), if it is both viscosity super- and
subsolution.
\end{definition}

We recall the function $u(t,x)$ is denoted by $u(t,x):=Y_{t}^{t,x}$, where $%
(Y_{s}^{t,x},Z_{s}^{t,x},A_{s}^{t,x})_{t\leq s\leq T}$ is the smallest
solution of BSDE (\ref{CBSDE}) constraint by $\Phi \geq 0$. And such
solution exists. Our first reault is following.

\begin{proposition}
\label{super-solu}For each $\alpha =1,2,\cdots $, $u(t,x)$ is a
discontinuous viscosity supersolution of $\partial _{t}u+F_{\alpha }=0$.
\end{proposition}

\noindent\textbf{Proof. } It is an application of Proposition
\ref{gmphi} and the fact that a $g$-supersolution relate to
viscosity supersolution. $\square $

Then we have

\begin{theorem}
The function $u\ $is a discontinuous viscosity solution of (\ref{PDECons}).
\end{theorem}

\noindent\textbf{Proof. }
From the above discussion, we know that for each $\alpha =1,2,\cdots $, $%
u_{\alpha }$, defined by $u_{\alpha }(t,x):=Y_{t}^{t,x,\alpha }$, is a
viscosity solution of $\partial _{t}u_{\alpha }+F_{\alpha }=0$. And $%
u_{\alpha }(t,x)\nearrow u(t,x)$, so $u(t,x)$ is lower semicontinuous, i.e. $%
u(t,x)=u_{*}(t,x)$.

We now prove that $u$ is a subsolution of (\ref{PDECons}). Let $(t,x)$ be a
point such that $\Phi (t,x,u^{*}(t,x),\sigma ^{T}q)>0$, and $(p,q,X)\in
\mathcal{P}^{2,+}u^{*}(t,x)$.

By Lemma 6.1 in \cite{CIL}, there exist sequences
\[
\alpha _{j}\rightarrow \infty ,(t_{j},x_{j})\rightarrow
(t,x),(p_{j},q_{j},X_{j})\in \mathcal{P}^{2,+}u_{\alpha _{j}}(t,x),
\]
such that
\[
(u_{\alpha _{j}}(t_{j},x_{j}),p_{j},q_{j},X_{j})\rightarrow
(u^{*}(t,x),p,q,X).
\]
While for each $j$,
\begin{eqnarray*}
&&-p_{j}-F_{\alpha }(t_{j},x_{j},u_{n_{j}}(t_{j},x_{j}),q_{j},X_{j}) \\
&=&-p_{j}-[F_{0}(t_{j},x_{j},u_{n_{j}}(t_{j},x_{j}),q_{j},X_{j})+\alpha \Phi
^{-}(t_{j},x_{j},u_{\alpha _{j}}(t_{j},x_{j}),\sigma ^{T}q_{j})]\leq 0.
\end{eqnarray*}
From the assumption that $\Phi (t,x,u^{*}(t,x),\sigma ^{T}q)>0$, continuity
assumption of $\Phi $ and convergence of $u_{\alpha }$, it follows for $j$
large enough, $\Phi ^{-}(t_{j},x_{j},u_{\alpha _{j}}(t_{j},x_{j}),\sigma
^{T}q_{j})>0$. Hence taking limit in the above inequality, we get
\[
-p-F_{0}(t,x,u^{*}(t,x),q,X)\leq 0,
\]
we prove that $u$ is viscosity subsolution of (\ref{PDECons}).

Then we conclude by proving that $u$ is a viscosity supersolution of (\ref
{PDECons}). Let $(t,x)\in [0,T]\times \mathbf{R}^{n}$, and $(p,q,X)\in
\mathcal{P}^{2,-}u_{*}(t,x)$, by proposition \ref{super-solu}, we know that $%
u$ is a discontinuous viscosity supersolution of $\partial _{t}u+F_{\alpha
}=0$, for each $\alpha \geq 0$, i.e.
\begin{eqnarray*}
&&-p-F_{\alpha }(t,x,u_{*}(t,x),q,X) \\
&=&-p-F_{0}(t,x,u_{*}(t,x),q,X)-\alpha \Phi ^{-}(t,x,u_{*},\sigma
^{T}(t,x)q)\geq 0.
\end{eqnarray*}
By the arbitrary of $\alpha $, we have
\[
-p-F_{0}(t,x,u_{*}(t,x),q,X)\geq 0\mbox{ and }\Phi
^{-}(t,x,u_{*},\sigma ^{T}(t,x)q)=0,
\]
i.e. $\Phi (t,x,u_{*}(t,x),\sigma ^{T}q)\geq 0$. $\square $

Then we consider the uniqueness of the solution. First, we have a
characterization property of $u(t,x)$.

\begin{proposition}
The function $u(t,x)$ is the smallest viscosity supersolution of (\ref
{PDECons}).
\end{proposition}

\noindent\textbf{Proof. }
Consider another viscosity supersolution of (\ref{PDECons}) denoted by $%
\overline{u}(t,x)$. By the definition, we have for each $(t,x)\in
(0,T)\times \mathbf{R}^{n}$, for any $(p,q,X)\in \mathcal{P}^{2,-}\overline{u%
}_{*}(t,x)$, then
\[
\Phi (x,\overline{u}_{*},\sigma ^{T}(x)q)\geq 0\mbox{ and }-p-F_{0}(t,x,%
\overline{u}_{*},q,X)\geq 0.
\]
So for $\alpha =1,2,\cdots $
\[
-p-F_{0}(t,x,\overline{u}_{*}(t,x),q,X)-\alpha \Phi ^{-}(t,x,\overline{u}%
_{*},\sigma ^{T}(t,x)q)\geq 0,
\]
which follows that $\overline{u}_{*}(t,x)$ is also a viscosity supersolution
of (\ref{eq4.9}). While $u_{\alpha }(t,x)$ is a viscosity solution of (\ref
{eq4.9}), then
\[
u_{\alpha }(t,x)\leq \overline{u}_{*}(t,x)\leq \overline{u}(t,x).
\]
By the limit property in (\ref{eq4.8}), we have $\overline{u}(t,x)\geq
u(t,x) $. And the result follows.

For the uniqueness of viscosity solution, we have following result.

\begin{theorem}
Under assumptions (\ref{Lip}) and (\ref{Constraint}), we assume that for
each $R>0$, there exists a function $m_{R}:\mathbf{R}_{+}\rightarrow \mathbf{%
R}_{+}$, such that $m_{R}(0)=0$ and
\[
\left| g(t,x_{1},y,p)-g(t,x_{2},y,p)\right| \leq m_{R}(\left|
x_{1}-x_{2}\right| (1+\left| p\right| ),
\]
for all $t\in [0,T]$, $\left| x_{1}\right| ,\left| x_{2}\right| \leq R$, $%
\left| y\right| \leq R$, $p\in \mathbf{R}^{d}$. And $\Phi $ is strictly
increasing in $y$ for each $(t,x,z)$. Then the constraint PDE \ref{PDECons}
has at most one locally bounded viscosity solution.
\end{theorem}

\noindent\textbf{Proof. } The proof is done by the same techniques
in theorem 8.6 in \cite{EKPPQ}, so we omit it. $\square $

\begin{remark}
The constraint satisfies assumptions in this theorem, if $\Phi
(t,x,y,z)=y-h(t,x)$, here $h(t,x)$ may be a discontinuous function with
certain integral condition. In fact such constraint introduces a reflected
BSDE with a discontinuous barrier $h(s,X_{s}^{t,x})$, c.f. \cite{PX2005}.
Another example is solution $y$ reflected on function of $z$, i.e. $\Phi
(t,x,y,z)=y-\varphi (t,x)$, where $\varphi $ is a Lipschitz function on $z$.
\end{remark}

\end{document}